\documentclass[final]{elsart}
\usepackage{latexsym,amsfonts,amssymb,amsbsy}
\usepackage{citesort}

\begin{document}
\begin{frontmatter}
\title{Series Prediction based on Algebraic Approximants}
\author[chemie,suc]{H. H. H. Homeier\thanksref{corresponding}}
\address[chemie]{Institut f\"ur Physikalische und Theoretische Chemie,
Universit\"at Regensburg, 93040 Regenburg, Germany}
\address[suc]{science $+$ computing ag, Ingolst\"adter Str. 22,
80807 M\"unchen, Germany}
\thanks[corresponding]{Email: Herbert.Homeier@na-net.ornl.gov, Phone:
+49-171-6290224, Fax: +49-89-356386737}
\journal{ISRN Applied Mathematics}


\maketitle

\begin{abstract} 
It is described how the Hermite-Pad\'e polynomials corresponding to
an algebraic approximant for a power series may be used to predict
coefficients of the power series that have not been used to compute the 
Hermite-Pad\'e polynomials. A recursive algorithm is derived and some numerical examples are given.
\end{abstract}
\end{frontmatter}

\noindent

\section{Introduction}
Using sequence transformation and extrapolation algorithms for the prediction of further sequence elements
from a finite number of known sequence elements is a topic of growing importance in applied mathematics. For a short introduction
see the book of Brezinski and Redivo Zaglia \cite[Sec. 6.8]{BrezinskiRedivoZaglia91}.
We mention theoretical work on prediction properties of Pad\'e approximants and related algorithms like the epsilon algorithm, and the iterated Aitken and Theta algorithms
\cite{Gilewicz73,PrevostVekemans99,Brezinski85b,Weniger00ppo}, Levin-type sequence transformations \cite{SidiLevin83,RoyBhattacharya06}, the E algorithm \cite{Brezinski85b,Vekemans97},  
and applications on
perturbation series of physical problems \cite{JentschuraBecherWenigerSoff00,RoyBhattacharya06}.

Here, we will concentrate on a different class of approximants, namely, the algebraic approximants. For a general introduction to these approximants and the related
Hermite-Pad\'e polynomials see
\cite{BakerGravesMorris96}. Programs for these approximants are available \cite{FeilHomeier04pft}. We summarize those properties that are important for the following.

Consider a function $f$ of complex variable $z$ with a known (formal)
power series
\begin{equation}
f(z)=\sum_{j=0}^{\infty} f_j \, z^j \>.
\end{equation}
The Hermite-Pad\'e polynomials (HPPs) corresponding to a certain algebraic
approximant are $N+1$
polynomials $P_n(z)$ with degree $p_n=deg(P_n)$, $n=0..N$  such that the
order condition
\begin{equation}\label{eqOrder}
\sum_{n=0}^N P_n(z) f(z)^n = O(z^M)
\end{equation}
holds for small $z$. Since one of the coefficients of the polynomials can be normalized
to unity, the order condition (\ref{eqOrder}) gives rise to a 
system of $M$ linear
equations for $N+\sum_{n=0}^N p_n$ unknown polynomial coefficients. Thus, the coefficient of $z^m$ of the Taylor expansion at $z=0$ of  the left hand side of Eq.~(\ref{eqOrder}) must be zero for $m=0,\dots,M-1$.
In order to have exactly as many equations as unknowns, we choose
\begin{equation}\label{eqNumpoints}
M=N+\sum_{n=0}^N p_n
\end{equation}
and assume that the linear system (\ref{eqOrder}) has a solution.
Then, the HPPs $P_n(z)$ are uniquely defined upon specifying the
normalization. The algebraic approximant under consideration then is that
pointwise solution $a(z)$ of the algebraic equation
\begin{equation}\label{eqAA}
P_0(z)+ \sum_{n=1}^N P_n(z) a(z)^n = 0
\end{equation}
for which the Taylor series of $a(z)$ coincides with the given power
series at least up to order $z^{M-1}$.

We note that for $N=1$, the algebraic approximants are nothing but the
well-known Pad\'e approximants.

Although we assumed that the power series of $f$ is known, quite often in
practice, only a finite number of coefficients \emph{is} really known.
These coefficients then may be used to compute the Hermite-Pad\'e
polynomials and the algebraic approximant under consideration.

We note that the higher coefficients of the Taylor series of $a(z)$ may be
considered as predictions for the higher coefficients of the power
series. The latter are also of interest in applications.

The question then arises how to compute the Taylor series of $a(z)$. If
it is possible to solve the equation (\ref{eqAA}) explicitly, i.e. for
$N\le 4$, a computer
algebra system may be used to do the job. But even then, a recursive
algorithm for the computation of the coefficients of the Taylor series
would be preferable in order to reduce computational efforts.

In the following section, such a recursive algorithm is obtained. In a
further section, we will present numerical examples.

\section{The recursive algorithm}

We consider the HPPs
\begin{equation}
P_n(z)=\sum_{j=0}^{p_n} p_{n,j} \, z^j
\end{equation}
as known. Putting
\begin{equation}
a(z) = \sum_{k=0}^{\infty} a_k \, z^k
\end{equation}
we obtain from Eq.\ (\ref{eqAA})
\begin{equation}\label{eqAAexp}
\sum_{j=0}^{p_0} p_{0,j} z^j + \sum_{n=1}^N \sum_{j=0}^{p_n} p_{n,j} z^j \sum_{k_1=0}^{\infty} \dots
\sum_{k_n=0}^{\infty} z^{k_1+\dots+k_n} \prod_{m=1}^{n} a_{k_m} = 0
\end{equation}
whence, by equating the coefficient of $z^J$ to zero, we obtain an
infinite set of equations. Due to Eq.\ (\ref{eqOrder}), all the equations
for $J<M$ are satisfied exactly for $a_j=f_j$, $j=0,\dots,M-1$. 

As a first step, we compute $a_M$.
We note that $M>p_0$. Hence, the coefficient of $z^M$ does not involve
any terms with $p_{0,j}$. For this coefficient $R_M$, we only need to
consider terms in Eq.\ (\ref{eqAAexp}) such that 
$M=j+k_1+\dots+k_n$ and we obtain $R_M=0$ for
\begin{equation}
R_M= \sum_{n=1}^N \sum_{j+k_1+\dots+k_n=M} p_{n,j} \prod_{m=1}^{n} a_{k_m}
\end{equation}
The only terms on the RHS involving $a_M$ are obtained if exactly one of
the $k_m$ is equal to $M$, i.e., we have $k_m=M$, $j=0$ and $k_j=0$ for
$j\ne m$. Thus, we may rewrite all these terms as $a_M C$ where
\begin{equation}\label{C}
C=\sum_{n=1}^N n\,p_{n,0} f_0^{n-1} 
\end{equation}
and note that the rest $D_M=R_M-a_M C$ is independent of $a_M$.
Recalling $R_M=0$, we obtain
\begin{equation}
a_M= - D_M / C
\end{equation}
Proceeding analogously for $J>M$, only terms with $J=j+k_1+\dots+k_n$
need to be considered. Hence, $R_J=0$ for
\begin{equation}
R_J= \sum_{n=1}^N \sum_{j+k_1+\dots+k_n=J} p_{n,j} \prod_{m=1}^{n} a_{k_m}
\end{equation}
Now, the only terms on the RHS involving $a_J$ are obtained if exactly
one of
the $k_m$ is equal to $J$, i.e., we have $k_m=J$, $j=0$ and $k_j=0$ for
$j\ne m$. Thus, we may rewrite all these terms as $a_J C$ where $C$ is
defined above. Proceeding as before, we put $D_J=R_J-a_J C$ and obtain
\begin{equation}\label{aJ}
a_J= - D_J / C
\end{equation}

An equivalent form of the recursive algorithm is obtained in the following way:

Consider for known $P_n$ and $a_0,\dots,a_{J-1}$ the expression
\begin{equation}
U_J=\left.\frac{d^J}{J! dz^J}\right\vert_{z=0} \sum_{n=1}^N P_n(z)\left(\sum_{j=0}^J a_j\,z^j\right)^n
\end{equation}
It is easy to see, that this expression is exactly equal to $R_J$, and hence, is linear in the unknown $a_J$. Thus, 
we may compute the quantities $D_J$ by substituting $a_J=0$ into $U_J$, which entails
\begin{equation}\label{DJ}
D_J=\left.\frac{d^J}{J! dz^J}\right\vert_{z=0} \sum_{n=1}^N P_n(z)\left(\sum_{j=0}^{J-1} a_j\,z^j\right)^n
\end{equation}

Thus, starting from $J=M$, one may compute all the $a_J$ consecutively by repeated use  of 
Eqs.\ (\ref{C}), (\ref{DJ}), and (\ref{aJ}).

This concludes the derivation of the recursive algorithm.

\section{Modes of application}

Basically, there are two modes of application.

a) one computes a sequence of HPPs and for the resulting algebraic approximants, one predicts a fixed number of sofar unused coefficients, e.g., only one new
coefficient. This mode is mainly for tests.

b) one computes from all available coefficients certain HPPs. For the best HPPs one computes a larger number of predictions for sofar unused coefficients.

In the following examples, we concentrate on mode b). Here, it is to be expected that the computed values have the larger errors the  higher coefficients are predicted.

\section{Examples}

The examples serve to introduce to the approach.
All numerical calculations in this section were done using Maple (Digits=16). 

\subsection*{Example 1}
As a first example, we consider $N=2$, $p_0=p_1=p_2=1$ and, hence, $M=5$.
Since $N=2$, we are dealing with a quadratic algebraic approximant.
Then, the recursive algorithm is started by
$a_j=f_j$, $j=0,\dots,4$. For $a_5$, we obtain
\begin{equation}
a_5= -
\frac
  {p_{1,1}f_4+p_{2,1}(2f_0f_4+2f_1f_3+f_2^2)+p_{2,0}(2f_1f_4+2f_2f_3)}
  {p_{1,0}+2p_{2,0}f_0}
\end{equation}
and for $J>5$, we obtain
\begin{equation}
a_J= -
\frac{p_{1,1}a_{J-1}+p_{2,1}\sum_{k=0}^{J-1}a_{J-k-1}a_k+p_{2,0}\sum_{k=1}^{J-1}
a_ka_{J-k}}{p_{1,0}+2p_{2,0}f_0}
\end{equation}

For 
\begin{equation}
f(z)= (2-3 z)^{1/2}+1/(5-z)
\label{ex1}
\end{equation}
the HPPs are determined  to be
\begin{eqnarray}
P_0(z)&{}={}&1.-1.544503593423590\, z \\
P_1(z)&{}={}&.1947992842134984+.06783822675080703\, z \\
P_2(z)&{}={}&-.5044536972622500-.01090573365920830\, z 
\end{eqnarray}

The results for the predicted coefficients given in Table \ref{tab1}.

\begin{table}[h]
\caption{The case of $N=2$, $p_0=p_1=p_2=1$ for Eq.\ (\ref{ex1}). Displayed are the coefficients of the Taylor series, the predicted coefficients and
absolute and relative errors of the predicted coefficients.}
\label{tab1}
\begin{tabular}{rrrrr}
$j$ & $f_j$ & $a_j$ & $\vert{f_j-a_j}\vert$ & rel. error (\%) \\ 
\hline
5 &     -.294 &     -.294 & .001 & .18 \\ 
6 &     -.330 &     -.332 & .001 & .38 \\ 
7 &     -.389 &     -.392 & .002 & .58 \\ 
8 &     -.475 &     -.478 & .004 & .76 \\ 
9 &     -.593 &     -.599 & .006 & .93 \\ 
10 &     -.756 &     -.765 & .008 & 1.10 \\ 
\end{tabular}
\end{table}

\subsection*{Example 2}
As a second example, we consider again $N=2$, $p_0=p_1=p_2=1$, and $M=5$, but now the function
\begin{equation}\label{ex2}
f(z)=17(1-2z)^{-1/3}+z/(2-z)
\end{equation}
with the HPPs
\begin{eqnarray}
P_0(z)&{}={}&-49.52369318166839-6.946105600281359\, z \\
P_1(z)&{}={}&1.+1.695055482965655\, z \\
P_2(z)&{}={}&.1125387307324166-.2732915349762758\, z
\end{eqnarray}

The results for the predicted coefficients given in Table \ref{tab2}.
\begin{table}[h]
\caption{The case of $N=2$, $p_0=p_1=p_2=1$ for Eq.\ (\ref{ex2}). Displayed are the coefficients of the Taylor series, the predicted coeffizients and
absolute and relative errors of the predicted coefficients.}
\label{tab2}
\begin{tabular}{rrrrr}
$j$ & $f_j$ & $a_j$ & $\vert{f_j-a_j}\vert$ & rel. error (\%) \\ 
\hline
5 &    67.938 &    68.212 & .274 & .40 \\ 
6 &   120.739 &   122.291 & 1.552 & 1.29 \\ 
7 &   218.459 &   224.194 & 5.735 & 2.62 \\ 
8 &   400.498 &   418.053 & 17.555 & 4.38 \\ 
9 &   741.657 &   790.063 & 48.406 & 6.53 \\ 
10 &  1384.425 &  1509.437 & 125.012 & 9.03 \\
\end{tabular}
\end{table}

\subsection*{Example 3}
As a final example,
we consider the case $N=p_0=p_1=p_2=2$, whence $M=8$, and the function
\begin{equation}\label{ex3}
f(z)=\exp(z)\,(2-3\,z)^{-1/3}+1/(5-z)
\end{equation}
The corresponding HPPs are
\begin{eqnarray}
P_0(z)&{}={}&1.-1.027576803009053\,z+.02070967420422950\,z^2 \\
P_1(z)&{}={}&2.617867885747464-.6563757889994458\,z \nonumber \\
      &     &{}-3.118191126500581\,z^2 \\
P_2(z)&{}={}&-3.647182626894738+7.471780741166546\,z \nonumber \\
      &     &{}-3.356878399103086\,z^2
\end{eqnarray}
The results for the predicted coefficients are displayed in Table \ref{tab3}.

\begin{table}[h]
\caption{The case of $N=2$, $p_0=p_1=p_2=2$ for Eq.\ (\ref{ex3}). Displayed are the coefficients of the Taylor series, the predicted coeffizients and
absolute and relative errors of the predicted coefficients.}
\label{tab3}
\begin{tabular}{rrrrr}
$j$ & $f_j$ & $a_j$ & $\vert{f_j-a_j}\vert$ & rel. error (\%) \\ 
\hline
8 &  3.888956 &  3.878509 &   .010447 &       .27 \\ 
9 &  5.356681 &  5.301047 &   .055634 &      1.04 \\ 
10 &  7.451679 &  7.275227 &   .176452 &      2.37 \\ 
11 & 10.447061 & 10.006950 &   .440111 &      4.21 \\ 
12 & 14.739132 & 13.781978 &   .957155 &      6.49 \\ 
13 & 20.903268 & 18.995972 &  1.907297 &      9.12 \\ 
\end{tabular}
\end{table}
\section{Conclusions}
It is seen that even rather low-order algebraic approximants, or HPPs, respectively, can lead to quite accurate predictions of the unknown coefficients of the power series,
especially for $f_M$, and the next few coefficients.

{\small
 \newcommand{\homeier}[1]{}

}


\begin{thebibliography}{10}

\bibitem{BakerGravesMorris96}
G.~A. {Baker, Jr.} and {P}. Graves-Morris.
\newblock {\em {P}ad{\'e} approximants}.
\newblock Cambridge U.P., Cambridge (GB), second edition, 1996.

\bibitem{Brezinski85b}
C.~Brezinski.
\newblock Prediction properties of some extrapolation methods.
\newblock {\em Appl. Numer. Math.}, 1:457 -- 462, 1985.

\bibitem{BrezinskiRedivoZaglia91}
C.~Brezinski and M.~{Redivo Zaglia}.
\newblock {\em Extrapolation methods. {T}heory and practice}.
\newblock North-Holland, Amsterdam, 1991.

\bibitem{FeilHomeier04pft}
T.~M. Feil and H.~H.~H. Homeier.
\newblock Programs for the approximation of real and imaginary single- and
  multi-valued functions by means of {Hermite-Pad\'{e}}-approximants.
\newblock {\em Comput. Phys. Commun.}, 158:124--135, 2004.
\newblock Computer Physics Communications Program Library, Catalogue number:
  ADSO.

\bibitem{Gilewicz73}
J.~Gilewicz.
\newblock Numerical detection of the best {Pad\'e} approximant and
  determination of the {Fourier} coefficients of insufficiently sampled
  functions.
\newblock In P.~R. Graves-Morris, editor, {\em {Pad\'e Approximants and their
  Applications}}, pages 99--103. Academic Press, New York, 1973.

\bibitem{JentschuraBecherWenigerSoff00}
U.~D. Jentschura, J.~Becher, E.~J. Weniger, and G.~Soff.
\newblock Resummation of {QED} perturbation series by sequence transformations
  and the prediction of perturbative coefficients.
\newblock {\em Phys. Rev. Lett.}, 85:2446--2449, 2000.

\bibitem{PrevostVekemans99}
M.~{Pr\'evost} and D.~Vekemans.
\newblock Partial {Pad\'e} prediction.
\newblock {\em Numer. Algo.}, 20:23--50, 1999.

\bibitem{RoyBhattacharya06}
Dhiranjan Roy and Ranjan Bhattacharya.
\newblock Prediction of unknown terms of a sequence and its application to some
  physical problems.
\newblock {\em Annals of Physics}, 321:1483--1523, 2006.

\bibitem{SidiLevin83}
A.~Sidi and D.~{L}evin.
\newblock Prediction properties of the {$t$}-transformation.
\newblock {\em SIAM J. Numer. Anal.}, 20:589--598, 1983.

\bibitem{Vekemans97}
D.~Vekemans.
\newblock Algorithm for the {E}-prediction.
\newblock {\em J. Comput. Appl. Math.}, 85:181--202, 1997.

\bibitem{Weniger00ppo}
E.~J. Weniger.
\newblock Prediction properties of {Aitken}'s iterated {$\Delta^2$} process, of
  {Wynn}'s epsilon algorithm, and of {Brezinski}'s iterated theta algorithm.
\newblock In C.~Brezinski, editor, {\em Numerical Analysis 2000 Vol. 2:
  Interpolation and Extrapolation}, pages 329 -- 356. Elsevier, Amsterdam,
  2000.

\end{thebibliography}
\end{document}